\documentclass[11pt,a4paper]{article}

\usepackage{amsmath, amssymb, enumerate, mathrsfs}
\usepackage[latin1]{inputenc} 
\usepackage[english,francais]{babel}
\usepackage{amsfonts}
\usepackage{amsthm}
\usepackage{relsize}
\usepackage{setspace}
\usepackage{geometry}
\usepackage{enumerate}
\usepackage{url}
\usepackage{xspace}
\usepackage{mathrsfs}
\usepackage{tocloft}

\usepackage{hyperref}

\usepackage{authblk}

\newtheorem{theorem}{Theorem}

\newtheorem{propo}[theorem]{Proposition}
\newtheorem{corollary}[theorem]{Corollary}
\newtheorem{e-definition}[theorem]{Definition\rm}

\newcommand{\lnorm}{\mathcal{LN}}
\newcommand{\lcent}{\mathcal{LC}}
\newcommand{\ld}{\mathcal{LD}}
\newcommand{\con}{\mathrm{con}}

\def\og{\leavevmode\raise.3ex\hbox{$\scriptscriptstyle\langle\!\langle$~}}
\def\fg{\leavevmode\raise.3ex\hbox{~$\!\scriptscriptstyle\,\rangle\!\rangle$}}

\title{Locally normal subgroups of simple locally compact groups}

\author[1]{Pierre-Emmanuel Caprace\thanks{F.R.S.-FNRS research associate, supported in part by the ERC (grant \#278469)}}
\author[2]{Colin D. Reid\thanks{Supported in part by ARC Discovery Project DP120100996}}
\author[2]{George A. Willis\thanks{Supported in part by ARC Discovery Project DP0984342}}

\affil[1]{Universit\'e catholique de Louvain, IRMP, Chemin du Cyclotron 2, bte L7.01.02, 1348 Louvain-la-Neuve, Belgique}
\affil[2]{Department of Mathematics, University of Newcastle, Callaghan, NSW 2308, Australia}

\date{March 27, 2013}

\begin{document}

\selectlanguage{english}
\maketitle

\begin{abstract}
We announce various results concerning the structure of compactly generated simple locally compact groups. We introduce a local invariant, called the \textbf{structure lattice}, which consists of commensurability classes of compact subgroups with open normaliser, and show that its properties reflect the global structure of the ambient group. 
\end{abstract}

\section{Introduction}
This note concerns the  class of compactly generated locally compact groups that are topologically simple and non-discrete. Since the connected component of the identity is always a closed normal subgroup in any topological group, the members that class are either connected or totally disconnected. The connected ones are known to coincide with the simple Lie groups, as a consequence of the solution to the Hilbert fifth problem. We shall therefore concentrate on the  totally disconnected ones. 
For the sake of brevity, we denote by $\mathscr S$ the class of compactly generated locally compact groups that are totally disconnected, topologically simple, and non-discrete. 

Our goal is to present some new tools to investigate the structure of  members of $\mathscr S$. Although these tools are of \emph{local} nature, i.e. depend only on arbitrarily small identity neighbourhoods, they interact with the global structure of the ambient group. The present study was inspired by earlier work due to J. Wilson \cite{Wilson} on just-infinite groups, and work by Barnea--Ershov--Weigel \cite{BEW} on abstract commensurators of profinite groups. For further considerations and detailed proofs, we refer to \cite{CRW}. 

\section{Locally normal subgroups}
The central  concept in our considerations is that of a \textbf{locally normal subgroup}, which is defined as a compact subgroup whose normaliser is open. Obvious examples are provided by the trivial subgroup, or by compact open subgroups. When $G$ is a simple Lie group or a simple algebraic group over a local field, every locally normal subgroup is of this form. However, in all other known examples of groups in $\mathscr S$, it has been observed, or is suspected, that there exist non-trivial locally normal subgroups that are not open. Our first result describes the algebraic structure of locally normal subgroups. 

\smallskip
\begin{theorem}\label{thm:primes}
Let $G \in \mathscr S$. Then:
\begin{enumerate}[(i)]
\item  There is a finite set of primes $\eta = \eta(G)$ such that every locally normal subgroup of $G$ is a virtually  pro-$\eta$ group. In particular the open pro-$\eta$ subgroups of $G$ form a basis of identity neighbourhoods. 
\item For every $p \in \eta$,  every locally normal subgroup $L \neq \{e\}$ has an infinite pro-$p$ subgroup.

\item If some non-trivial locally normal subgroup is virtually pro-soluble, then so is every compact open subgroup. 
\item The only virtually soluble locally normal subgroup is the trivial subgroup $\{e\}$.
\end{enumerate}
\end{theorem}

\section{The structure lattice}
We next consider the set $\lnorm(G)$  of all locally normal subgroups, modulo the equivalence relation defined by commensurability. We endow $\lnorm(G)$ with the partial order $\leq$ induced by inclusion of locally normal subgroups. One verifies that $\lnorm(G)$ is a modular lattice, endowed with a canonical $G$-action by automorphisms, induced by the conjugation action. The \emph{join} and \emph{meet} operations in $\lnorm(G)$ are denoted by $\vee$ and  $\wedge$ respectively. We call $\lnorm(G)$ the \textbf{structure lattice} of $G$.  It possesses a global minimum $0$, corresponding to
the class of the trivial subgroup, and a global maximum $\infty$,
corresponding to the class of open compact subgroups.
By \cite[Th.~4.8]{BEW} (or as a consequence of Theorem~\ref{thm:primes}(iv) above), the only finite locally normal subgroup is the trivial one. It follows that $\lnorm(G) = \{0, \infty\}$ if and only if all compact open subgroups of $G$ are hereditarily just-infinite. This is automatically the case when $G$ is a simple algebraic group over a local field. 

Our next result provides information on the set of fixed points of $G$ in the structure lattice, which we denote by $\lnorm(G)^G$. 

\smallskip
\begin{theorem}\label{thm:fp}
Let $G \in \mathscr S$. Then:
\begin{enumerate}[(i)]
\item There is a unique element $\mu \in \lnorm(G)$ such that  $ \lnorm(G)^G \setminus \{0\} = \{\alpha \in \lnorm(G) \; | \; \alpha \geq \mu\}$.

\item $\mu = \infty$ if and only if for each non-trivial locally normal subgroup $K$, there are finitely many conjugates $K_1, \dots, K_n$ of $K$ such that the product $K_1. \dots . K_n$ is an identity neighbourhood.

\item If $\mu= \infty$, then every compact subgroup of $G$ commensurated by $G$ is either finite or open. 


\item If $G$ is abstractly simple, then $\mu = \infty$. 
\end{enumerate}
\end{theorem}
\smallskip

All known examples of groups in $\mathscr S$ have been proved to be abstractly simple, and it is tempting to believe that this is always the case.  This would imply that the only fixed points of $G$ in the structure lattice are the trivial ones, and that the only infinite commensurated compact subgroups  are open. 

The proof of Theorem~\ref{thm:fp}(i) relies on recent work by Nikolov--Segal \cite{NS}, which in turn depends on the classification of the finite simple groups.


\section{The centraliser lattice}
We next introduce an operator $\perp$ on $\lnorm(G)$ related to centralisers. Let $G \in \mathscr S$ and let $\alpha \in \lnorm(G)$. We may find a representative $K$ of $\alpha$ that is a closed normal subgroup of some compact open subgroup $U$ of $G$. The centraliser $\mathrm C_U(K)$ is then also a closed normal subgroup of $U$, and is thus itself a locally normal subgroup of $G$. Moreover, it can be deduced from   Theorem~\ref{thm:primes}(iv) that the group $\mathrm C_U(K) $ does not depend on the choice of the representative $K$ of $\alpha$. Therefore the commensurability class of $\mathrm C_U(K)$ depends only on $\alpha$; we denote it by $\alpha^\perp$. Clearly $0^\perp = \infty$ and $\infty^\perp = 0$. We set 
$\lcent(G) = \{\alpha^\perp \; | \; \alpha \in \lnorm(G)\}$ and call it the \textbf{centraliser lattice}. It follows from Theorem~\ref{thm:primes}(iv) that for each $\alpha \in \lnorm(G)$, we have $\alpha \wedge \alpha^\perp = 0$. It is however not true in general that $\alpha \vee \alpha^\perp = \infty$. In other words, it is not necessarily true that $\alpha$ and $\alpha^\perp$ admit two representatives whose product is an open subgroup of $G$. To remedy this fact, we introduce an abstract operator $\vee'$ defined by $\alpha \vee' \beta = (\alpha^\perp \wedge \beta^\perp)^\perp$. The interest of that definition is revealed by the following. 

\smallskip
\begin{theorem}\label{thm:Boole}
Let $G \in \mathscr S$. Then  the centraliser lattice $\lcent(G)$ endowed with the operators $\vee'$, $\wedge$ and $\perp$, is a Boolean lattice on which $G$ acts by automorphisms. 
\end{theorem}
\smallskip

Just as the structure lattice is a local object, so is the centraliser lattice: it can be entirely reconstructed from any identity neighbourhood in $G$. 
By the Stone representation theorem, any Boolean lattice is canonically isomorphic to the lattice of clopen sets of a compact totally disconnected space, called the \textbf{Stone space} of that lattice. In particular the automorphism group of the lattice is canonically isomorphic to the homeomorphism group of the Stone space. We shall denote the Stone space of $\lcent(G)$ by $\Omega(G)$. Notice that $\Omega(G)$ reduces to a singleton if and only if $\lcent(G) = \{0, \infty\}$. Our next result shows that, as soon as this is not the case, the dynamics of the $G$-action on $\Omega(G)$ is quite rich. 

\smallskip
\begin{theorem}\label{thm:dyn}
Let $G \in \mathscr S$. Assume that  $\lcent(G) \neq \{0, \infty\}$ and set $\Omega = \Omega(G)$. Then:
\begin{enumerate}[(i)]
\item $\Omega$ has no isolated point. 

\item The $G$-action on $\Omega$ is continuous and faithful. 

\item The $G$-action on $\Omega$ is \textbf{minimal}, i.e. every $G$-orbit   is dense. 

\item The $G$-action on $\Omega$ is \textbf{strongly proximal}, i.e. the closure of every $G$-orbit in the space of probability measures on $\Omega$, endowed with the weak-$*$ topology, contains a Dirac mass.

\item Every point of $\Omega$ has a {compressible} neighbourhood.

\end{enumerate}
\end{theorem}
\smallskip

A subset $V$ of $\Omega$ is called \textbf{compressible} if there is a sequence $(g_n)$ in $G$ such that $g_n.V$ converges to a singleton in the space of closed subsets of $\Omega$.   The existence of a non-trivial strongly proximal actions is clearly a strong obstruction to amenability. In fact, using Furstenberg's boundary theory \cite{Fur}, we deduce the following:

\smallskip
\begin{corollary}\label{cor:amen}
Let  $G \in \mathscr S$. Any closed cocompact amenable subgroup of $G$ fixes a point in $\Omega$. In particular, if  $\lcent(G) \neq \{0, \infty\}$, then $G$ is not amenable, and if $G$ contains a closed cocompact amenable subgroup, then the $G$-action on $\Omega$ is transitive.
\end{corollary}
\smallskip

This result contrasts with the recent groundbreaking work by Juschenko--Monod \cite{JM}, who obtained the first examples of finitely generated infinite simple groups that are amenable. It is not known whether non-discrete analogues of such groups exist or, more precisely, whether a group $G \in \mathscr S$ can be amenable. Corollary~\ref{cor:amen} provides some evidence that this might not be the case. 

Let us also point out the following fact about the topology of $G$ in case the centraliser lattice is non-trivial:

\smallskip
\begin{theorem}
Let $G \in \mathscr S$. If  $\lcent(G) \neq \{0, \infty\}$, then the topology of $G$ is the unique $\sigma$-compact locally compact group topology on $G$. In particular every automorphism of $G$ is continuous.  
\end{theorem}
\smallskip

The condition that the centraliser lattice is not $\{0, \infty\}$ is equivalent to the condition that some locally normal subgroup of $G$ splits non-trivially as the direct product of two locally normal subgroups; in fact the elements of $\lcent(G)$ account for direct factors of locally normal subgroups in such direct decompositions. Theorem~\ref{thm:dyn} can also be used to derive another algebraic  characterisation of those groups whose centraliser lattice is non-trivial:

\smallskip
\begin{propo}
Let $G \in \mathscr S$. Then $\lcent(G) \neq \{0, \infty\}$ if and only if $G$ has a closed subgroup isomorphic to the unrestricted wreath product $H = (\prod_{\mathbf Z} K) \rtimes \mathbf Z$, where $K$ is a non-trivial locally normal subgroup of $G$. 
\end{propo}
\smallskip

The `if' part is clear, since any two distinct conjugates of $K$ in $H$ provide two commuting locally normal subgroups of $G$, thereby ensuring the non-triviality of the centraliser lattice. The converse part of the proposition has the following noteworthy consequence: 

\smallskip
\begin{corollary}\label{cor:con}
Let $G \in \mathscr S$. If $\lcent(G) \neq \{0, \infty\}$, then the contraction group of some element  of $G$ is not closed. 
\end{corollary}
\smallskip

Recall that the \textbf{contraction group} of $g$ is defined by 
$$\con(g) = \{ x \in G \; | \; \lim_{n \to \infty} g^n x g^{-n} = e\}.$$ 
In simple Lie groups or in simple algebraic groups over local fields, the contraction group of every element is known to coincide with the unipotent radical of some parabolic subgroup, and is thus always closed, while for many non-linear examples of groups in $\mathscr S$, the existence of non-closed contraction group has been observed by means of a case-by-case analysis. Corollary~\ref{cor:con} provides strong evidence that this phenomenon should have a more conceptual explanation.


Another consequence of Theorem~\ref{thm:dyn} concerns abstract simplicity. We pointed out above that Theorem~\ref{thm:fp} can be strengthened  if some non-trivial locally normal subgroup of $G \in \mathscr S$ is finitely generated. This is also the case  if the centraliser lattice is non-trivial:

\smallskip
\begin{corollary}\label{cor:fg}
Let $G \in \mathscr S$. Assume that some non-trivial locally normal subgroup is finitely generated, or that $\lcent(G) \neq \{0, \infty\}$. Then  $G$ is abstractly simple if and only if $\lnorm(G)^G = \{0, \infty\}$. 
\end{corollary}

\section{The decomposition lattice}

There is a situation where the conclusion of Corollary~\ref{cor:fg} can be further improved. In order to describe it, we introduce yet another lattice related to $\lnorm(G)$, defined as follows: 
$\ld(G) = \{\alpha \in \lnorm(G) \; | \;  \alpha \vee \alpha^\perp = \infty\}.$ 
One verifies that  for $\alpha \in \ld(G)$ one has $(\alpha^\perp)^\perp = \alpha$; in particular $\ld(G)$ is contained in $\lcent(G)$. The lattice $\ld(G)$ is called the \textbf{decomposition lattice}. By definition, the restriction to $\ld(G)$ of the operator $\vee'$ introduced above coincides with $\vee$. Therefore, it follows from Theorem~\ref{thm:Boole}  that $\ld(G)$, endowed with the operators $\vee, \wedge$ and $\perp$, is itself a Boolean lattice. Its elements account for direct factors of compact open subgroups of $G$. In particular $\ld(G) = \{0, \infty\}$ if and only if no open subgroup of $G$ admits a decomposition as a direct product of two non-trivial closed factors. 

\smallskip
\begin{theorem}\label{thm:ld}
Let $G \in \mathscr S$ be such that some open subgroup of $G$ admits a decomposition as a direct product of two non-trivial closed factors. Then $G$ is abstractly simple. 
\end{theorem}
 

\section{Five types of simple groups}

 We finally collect some of the information about the structure lattice and its sub-lattices introduced above, in order to partition the class $\mathscr S$ into five distinct types, defined according to the properties of the structure lattice. We use the term \textbf{atom} in a lattice to qualify  a non-zero element which is not minorized by any other non-zero element. Accordingly, a lattice is called \textbf{non-atomic} if it does not contain any atom; in the case of a Boolean lattice, this is equivalent to the statement that the associated Stone space does not have isolated points. 

\smallskip
\begin{theorem}\label{thm:types}
Let $G \in \mathscr S$. Then exactly one of the following holds. 
\begin{enumerate}[(a)]
\item $\lnorm(G) = \{0, \infty\}$, and every compact open subgroup of $G$ is hereditarily just-infinite.

\item $\lnorm(G)$ is infinite, non-atomic,  the $G$-action on $\lnorm(G)$ is faithful, and $\lcent(G) = \{0, \infty\}$.

\item $\lcent(G)$ is infinite, non-atomic, and $\ld(G) = \{0, \infty\}$. 

\item $\ld(G)$ is infinite and non-atomic. 

\item $\lnorm(G) \neq  \{0, \infty\}$,  $\lcent(G) = \{0, \infty\}$ and the $G$-action on $\lnorm(G)$ is trivial. 
\end{enumerate}
\end{theorem}
\smallskip

In case (a) (resp. (b), (c), (d), and (e)), we say that $G$ is of  \textbf{h.j.i. type} (resp. \textbf{faithful type}, \textbf{weakly branch type}, \textbf{branch type}, and  \textbf{atomic type}). In case $G$ is of atomic type, the element $\mu$ from Theorem~\ref{thm:fp} is necessarily the unique atom in $\lnorm(G)$. Moreover Theorem~\ref{thm:fp} implies that a group of atomic type cannot be abstractly simple. We do not know whether groups of atomic type exist. On the other hand, the other four classes are  non-empty: groups of faithful type may be found among complete Kac--Moody groups over finite fields, groups of weakly branch type among automorphism groups of right-angled buildings, and groups of branch type among automorphism groups of trees. We point out a last consequence of our results on the algebraic structure of locally normal subgroups, which supplements Theorem~\ref{thm:primes}:

\smallskip
\begin{corollary}
Let $G \in \mathscr S$. Any hereditarily just-infinite locally normal subgroup is commensurated by $G$. Moreover, if $G$ contains such a subgroup, then $G$ is either of h.j.i. type or of atomic type. 
\end{corollary}

\end{document}